\documentclass[12pt]{article}
\usepackage{amsmath, amsfonts, amsthm,amssymb, color, hyperref, enumerate, extarrows, cite, subfiles, csquotes, enumitem, commath, mathtools,eucal, }

 \usepackage[titletoc]{appendix}

\newtheorem{Lemma}{Lemma}[section]
\newtheorem{Theorem}[Lemma]{Theorem}
\newtheorem{Proposition}[Lemma]{Proposition}
\newtheorem{Corollary}[Lemma]{Corollary}

\newtheorem{remark}[Lemma]{Remark}
\newtheorem{definition}[Lemma]{Definition}
\newtheorem{example}[Lemma]{Example}

\newtheorem{Fact}[Lemma]{Fact}

\hypersetup{
	colorlinks   = true,
	citecolor    = magenta}

\usepackage[titletoc]{appendix}

\def\bt{\begin{Theorem}}
\def\et{\end{Theorem}}
\def\bl{\begin{Lemma}}
\def\el{\end{Lemma}}
\def\bp{\begin{Proposition}}
\def\ep{\end{Proposition}}
\def\bcor{\begin{Corollary}}
\def\ecor{\end{Corollary}}
\def\bpf{\begin{proof}}
\def\epf{\end{proof}}

\def\brem{\begin{remark}}
\def\erem{\end{remark}}

\def\bedef{\begin{definition}\rm }
\def\endef{\end{definition}}

\def\beg{\begin{example}}
\def\eeg{\end{example}}

\def\bef{\begin{Fact}}
\def\eef{\end{Fact}}

\def\bc{\begin{center}}
\def\ec{\end{center}}

\def\beq{\begin{equation}}
\def\eeq{\end{equation}}
\def\beqarray{\begin{eqnarray*}}
\def\eeqarray{\end{eqnarray*}}
\def\<{\leftangle}
\def\>{\rightangle}
\def\({\left(}
\def\){\right)}

\def\<{\langle}
\def\>{\rangle}

\def\r{\rho}
\def\a{\alpha}

\def\p{\pi}
\def\d{\delta}

\def\k{\kappa}

\def\t{\tau}
\def\r{\rho}

\def\e{\varepsilon}

\def\O{\Omega}

\def\z{\zeta}

\def\w.r.t.{with respect to}
\def\R{{\mathbb{R}}}
\def\N{{\mathbb{N}}}
\def\Z{{\mathbb{Z}}}

\def\P{{\mathbb{P}}}
\def\C{{\mathbb{C}}}

\def\A{{\mathcal{A}}}
\def\L{{\mathcal{L}}}

\def\bq{\begin{quote}}
\def\eq{\end{quote}}
\renewcommand{\Re}{\mathop{\textrm{Re}}}

\def\bit{\begin{itemize}}
\def\eit{\end{itemize}}

\def\ben{\begin{enumerate}}
\def\een{\end{enumerate}}

\begin{document}

\title{\vspace{-1.2cm} \bf $p$-Skwarczy\'nski distance\rm}

\author{Shreedhar Bhat}
\date{}

\maketitle

\begin{abstract}

We introduce a new distance on a domain $\O \subset \C^n$ using the `minimizer' functions on $\A^p(\O)$. We discuss its invariance, completeness and other aspects related to it.

\end{abstract}

\begin{NoHyper}
\renewcommand{\thefootnote}{\fnsymbol{footnote}}
\footnotetext{\hspace*{-7mm} 
\begin{tabular}{@{}r@{}p{16.5cm}@{}}
& Keywords. Bergman kernel, Skwarczy\'nski distance, $p$-Bergman kernel.\\
& Mathematics Subject Classification. Primary 32F45; Secondary 32A36\end{tabular}}
\end{NoHyper}

\section{Introduction}
The study of invariant distances on complex domains has been a fundamental topic in complex analysis for decades. The Skwarczy\'nski distance, introduced in 1980 \cite{skwarczynski1980biholomorphic}, has been extensively studied due to its invariance with respect to biholomorphic functions and its close relationship to the Bergman metric \cite{mazur1985invariant}, \cite{pflug1982various}, \cite{skwarczynski1985general}, \cite{JarnickiPflug+2013}. Recently, a similar distance was defined by Krantz et al.\  in 2021, which utilized the Szegő kernel\cite{krantz2021invariant}. In this paper, we introduce a new distance on a bounded domain in $\C^n$ using the `minimizer' function on the $p$-Bergman space $\A^p(\O)$. We investigate the properties of this new distance, invariance, completeness, and potential applications of the $p$-Skwarczy\'nski distance.\\

For a bounded domain $\Omega\subset \C^n$, $1\leq p<\infty $, the $p$-Bergman space is defined as $$\A^p(\O)=\{f\in \mathcal{L}^p(\O):f \text{ is holomorphic in } \O\}$$
Define $$m_p(z_0)=\inf\{\norm{f}_p: f\in \A^p(\O), f(z_0)=1\}$$
From \cite[Proposition 2.4, 2.5]{chen2022p}, we know that there exists a unique such function which minimizes the above norm. Let the unique minimizer function be $m_p(\cdot,z_0)$. Define $$  \text{ the }p\text{-Bergman kernel }K_p(z_0)=m_p(z_0)^{-p}; $$ 
$$\text{ and the off-diagonal }p \text{-Begrman kernel }  K_p(\cdot,z_0)= K_p(z_0)m_p(\cdot,z_0)$$
 The above function mimics certain properties of the Bergman kernel ($K_2$) in $\A^p(\O)$ including a reproducing formula:
 \begin{align*}
     f(z)=m_p(z)^{-p}\int_\O \abs{m_p(w,z)}^{p-2} \overline{m_p(w,z)}f(w)dw \text{ for } f\in \A^p(\O)
 \end{align*}For more properties of the above function cf.\  \cite{chen2022p}, \cite{ning2016p}, \cite{chen2022some}, \cite{chen2023regularity}.\\



\section{\texorpdfstring{$p$}{p}-Skwarczy\'nski distance}\label{SecDefinition}

On $\A^p(\O) \setminus\{0\}$ define a relation $\sim$ by $f\sim g$ if and only if $f=c\cdot g$ for some non-zero complex number $c$. Let $\P(\A^p(\O))=[\A^p(\O)\setminus \{0\}]/\sim$ denote the projective space\\
Let $S_{\A_p}$ denote the unit sphere of $\A^p$. Equip $\P(\A^p)$ with the distance 
\begin{align*}
d([f],[g])&=\text{dist}([f]\cap S_{\A^p}, [g]\cap S_{\A^p}) 
\\ &=\inf_{t_1,t_2 \in [0,2\p]} \norm{ e^{i t_1}\frac{f}{\norm{f}}- e^{it_2}\frac{g}{\norm{g}}}_p\\
&=\norm{e^{it}\frac{f}{\norm{f}}-\frac{g}{\norm{g}}}_p \qquad \text{ for some } t \in \R. \end{align*}

\bl 
$(\P(\A^p(\O)) ,d)$ is a complete metric space.
\el 

\bpf 
Let $[f_n]$ be a Cauchy sequence and assume without loss of generality that $\norm{f_n}=1$ for every $n$.
We can find a subsequence  $[f_{n_k}]$ such that
$d([f_{n_k}],[f_{n_{k+1}}])\leq \frac{1}{2^k}$. 

Let $g_1 \in [f_{n_1}]$, $\norm{g_1}=1$ such that $dist (g_1,[f_{n_2}])\leq \frac{1}{2}$.

Let $g_2\in [f_{n_2}]$, $\norm{g_2}=1$ such that $\norm{g_1-g_2}= d([f_{n_1}],[f_{n_2}]) \leq 1/2$.

\hspace{1mm}\vdots 

Let $g_k \in [f_{n_k}]$, $\norm{g_k}=1$ such that $\norm{g_{{k-1}}-g_{k}}= d([f_{n_{k-1}}],[f_{n_k}])\leq \frac{1}{2^{k-1}}$.\\
Then for $r<s$, 
\begin{align*}
\norm{g_{r}-g_{s}} &\leq \norm{g_{r}-g_{r+1}}+\dots+\norm{g_{s-1}-g_{s} }\\
\hspace{-25mm}&\leq \frac{1}{2^{r}}+\dots +\frac{1}{2^{s-1}}=\left(\frac{1}{2}\right)^{r-1}-\left(\frac{1}{2}\right)^{s-1} \xrightarrow{} 0 \text{ as } r,s \rightarrow \infty.
\end{align*}
Thus $\{g_r\}$ is a Cauchy sequence in $\A^p(\O)$ and hence there exists a function $f\in \A^p(\O)$ , $\norm{f}=1$ such that $g_r \rightarrow f.$\\By definition 
$$d([f_{n_r}],[f])=d([g_r],[f]) \leq \norm{g_r-f} \rightarrow 0$$
Using the elementary property that if a subsequence of a Cauchy sequence converges, then the sequence converges,  $[f_n] \xrightarrow{d} [f]$, that is, every Cauchy sequence is a convergent sequence. 
\epf

\paragraph{}

Consider the map $$\t : \O\rightarrow \P(\A^p(\O)) $$ $$\hspace{25mm} z \longmapsto [K_p(\cdot,z)]=[m_p(\cdot,z)]$$ 
From \cite[Proposition 2.15]{chen2022p}, we know that $m_p(\cdot,z)\neq c \cdot m_p(\cdot,w)$ for any $c\in \C$, when $z\neq w \in \O$.\\
The above map is an injective map, so we can pull back the distance on $\t(\O)$ onto $\O$, that is,
\begin{align*}
\r_p(z,w)&=\text{dist}([K_p(\cdot,z)],[K_p(\cdot,w)])= \text{dist}([m_p(\cdot,z)],[m_p(\cdot,w)])\\ 
&=\inf_{\theta_1,\theta_2 \in \R }\norm{\frac{e^{i\theta_1} m_p(\cdot,z)}{m_p(z)}-\frac{e^{i\theta_2} m_p(\cdot,w)}{m_p(w)}}_p \\ & =\norm{\frac{e^{i\theta} m_p(\cdot,z)}{m_p(z)}-\frac{ m_p(\cdot,w)}{m_p(w)}}_p \quad \text{ for some } \theta \in [0,2\pi].
\end{align*}

\bedef 
The distance $\r_p$ defined above is called the \textbf{$p$-Skwarczy\'nski distance}. Thus 
\beq \label{SkwaDef}
\r_p(z,w)= \min_{t\in [0,2\pi]}\norm{e^{it}\frac{m_p(\cdot,z)}{m_p(z)}-\frac{m_p(\cdot,w)}{m_p(w)}}_p.
\eeq 
\endef

\section{Properties of the \texorpdfstring{$p$}{p}-Skwarczy\'nski distance} \label{SecProperties}
\subsection{Comparing the Euclidean topology and \texorpdfstring{$p$}{p}-Skwarczy\'nski topology}
\cite[Theorem III.8]{skwarczynski1980biholomorphic} states that the Euclidean topology and the Skwarczy\'nski topology agree for $p=2$. We will now investigate these topologies for $p \in [1,\infty)$   
\bl \label{EucImpSkwa}
Let $\{w_k\}$ be a sequence in $\O$ such that $w_k\xrightarrow{\text{Euclidean}}w$. Then $w_k$ converges to $w$ in the $p$-Skwarczy\'nski topology.
\el
\bpf
This follows directly from the following lemma \cite[Theorem 1]{brazis-lieb}.
\begin{displayquote}
\textbf{Brezis--Lieb Lemma: } For $p\in (0,\infty)$, if $f_k, f\in \L^p$ satisfy $f_k\xrightarrow{a.e.} f$, and $\norm{f_k}_p \rightarrow \norm{f}_p$, then $\norm{f_k-f}_p\rightarrow 0$.
\end{displayquote}
From \cite{chen2022p}, we know that $\frac{m_p(\cdot,w_k)}{m_p(w_k)} \xrightarrow{pointwise}\frac{m_p(\cdot,w)}{m_p(w)}$ and $\norm{\frac{m_p(\cdot,w_k)}{m_p(w_k)} }=\norm{\frac{m_p(\cdot,w)}{m_p(w)}}=1$. Thus $\r_p(w_k,w) \rightarrow 0$.
\epf 

\bcor 
The $p$-Skwarczy\'nski distance $\r_p$ is continuous, that is, $w_k \xrightarrow{Euclidean} w \implies \r_p(w_k,z)\rightarrow \r_p(w,z)$ for every $z$ in~$\Omega$.
\ecor 
\bpf
By triangle inequality, $\abs{\r_p(w_k,z)-\r_p(w,z)}\leq \r_p(w_k,w)\rightarrow 0$.
\epf

\bl \label{SkwaImpEuc}
If $p\geq 2$, and $\{w_k\}$ is a sequence in $\O$ such that $w_k\xrightarrow{\r_p}w$, then $w_k$ converges to $w$ in the Euclidean topology.
\el 
\bpf
Suppose $w_k\xrightarrow{\r_p} w$ in $\r_p$, equivalently, 
$$ \norm{\frac{e^{i\theta_k}m_p(\cdot,w_k)}{m_p(w_k)}-\frac{m_p(\cdot,w)}{m_p(w)}}_p \rightarrow 0 \qquad\text{ for some }\{\theta_k\}.$$ 

Then \begin{align*} 
h_k=e^{i\theta_k}\frac{m_p(\cdot,w_k)}{m_p(w_k)} \xrightarrow{normally}h=\frac{m_p(\cdot,w)}{m_p(w)}.\\ \end{align*}
For $p\geq 2$, $\abs{h_k}^{p-2}h_k$ converges normally to $\abs{h}^{p-2}h$. Observing that $C_c^\infty(\O)$ functions are dense in $\L^p(\O)$, we get that normal convergence implies weak convergence i.e.
\begin{align*}
    \< h_k \abs{h_k}^{p-2}, g\> \xrightarrow{}\< h\abs{h}^{p-2}, g\> \text{ for all } g\in \L^p(\O). \end{align*}
Since the domain $\O$ is bounded, $\{1, \pi_i\}\in \L^p(\O)$, where $\pi_i$ is the coordinate projection onto the $i^{th}$ coordinate, $\pi_i((w_1,\dots,w_n))=w_i$. Using the `reproducing formula' from \cite{chen2022p},
$$\lim_{m\rightarrow\infty}  {e^{-i\theta_k}} m_p(w_k) = \lim_{k\rightarrow\infty}  \int_\O \abs{\frac{m_p(\cdot,w_k)}{m_p(w_k)}}^{p-2}\overline{e^{i\theta_k}\frac{m_p(\cdot,w_k)}{m_p(w_k)}} = \int_\O \abs{\frac{m_p(\cdot,w)}{m_p(w)}}^{p-2}\overline{\frac{m_p(\cdot,w)}{m_p(w)}}=m_p(w).$$
Moreover, for $i=1,\dots,n$,
\begin{align*}\lim_{k\rightarrow\infty} \pi_i(w_k)  {e^{-i\theta_k}} m_p(w_k) &= \lim_{k\rightarrow\infty}  \int_\O \abs{\frac{m_p(\cdot,w_k)}{m_p(w_k)}}^{p-2}\overline{e^{i\theta_k}\frac{m_p(\cdot,w_k)}{m_p(w_k)}} \pi_i \\&= \int_\O \abs{\frac{m_p(\cdot,w)}{m_p(w)}}^{p-2}\overline{\frac{m_p(\cdot,w)}{m_p(w)}} \pi_i=m_p(w) \pi_i(w).\end{align*}
Thus $w_k\xrightarrow{Euclidean } w$.
\epf

\bl
If $p>2$, then the $p$-Skwarczy\'nski distance is locally $\frac{1}{p}$-H\"older continuous, that is,  
$$\forall z_0 \in\O, \hspace{2mm}\exists M,r>0 : \r_p(z,w)\leq M\abs{z-w}^{1/p} \qquad z,w\in B(z_0,r)\subset \O.$$
\el

\bpf
Let $z_0\in \O$ and $r>0$ such that $B(z_0,r)\subset\subset \O$. Let $z,w\in B(z_0,r)$.

We use the following inequality from \cite[Proposition 4.3 (1)]{chen2022p}:
\beq \label{4.3.1}
\abs{b-a}^p \leq 2^{p-1} \left[ \abs{b}^p+\abs{a}^p-\Re(\abs{b}^{p-2}\Bar{b}a+ \abs{a}^{p-2}\Bar{a}b)\right] \text{ when } p\geq 2.
\eeq 

Set $b=\frac{m_p(\z,z)}{m_p(z)}$, $a=\frac{m_p(\z,w)}{m_p(w)}$, and integrate over $\O$.

$$ 
\int_\O \abs{\frac{m_p(\z,z)}{m_p(z)}-\frac{m_p(\z,w)}{m_p(w)}}^p d\z$$ 
$$\leq 2^{p-1} \left[ 1+1-\Re\left[\int_\O\abs{\frac{m_p(\z,z)}{m_p(z)}}^{p-2} \overline{\frac{m_p(\z,z)}{m_p(z)}} \frac{m_p(\z,w)}{m_p(w)}d\z +\int_\O\abs{\frac{m_p(\z,w)}{m_p(w)}}^{p-2} \overline{\frac{m_p(\z,w)}{m_p(w)}} \frac{m_p(\z,z)}{m_p(z)} d\z\right] \right]
$$ 
By the reproducing formula,
$$\r_p(z,w)^p \leq 2^{p-1} \Re\left[ \frac{m_p(w)-m_p(z,w)m_p(z)}{m_p(w)}+\frac{m_p(z)-m_p(w,z)m_p(w)}{m_p(z)} \right]$$

Let $F(z,w)=\frac{m_p(w)-m_p(z,w)m_p(z)}{m_p(w)}$. Then $F(w,w)=0$, and from \cite{chen2022some}, we know that $m_p(\cdot), m_p(\cdot,w)$ are $C^1$ functions. Thus $F(\cdot,w)$ is a $C^1$ function, hence locally Lipschitz. Thus
$$\abs{F(z,w)-F(w,w)}=\abs{F(z,w)}\leq C_{z_0}\abs{z-w} \text{ for } z,w\in B(z_0,r) \subset\subset \O.$$
where $C_{z_0}=\max\left\{\abs{\frac{d}{dz}F(z,w)}: z,w\in \overline{B(z_0,r)}\right\}$.\\
Accordingly, $\r_p(z,w)^p \leq 2^{p-1}\left[\abs{F(z,w)}+\abs{F(w,z)}\right] \leq C' \abs{z-w}$, which implies the conclusion of the lemma.
\epf

\bl
If $p\in (1,2) $, then the  $p$-Skwarczy\'nski distance is locally $\frac{1}{2}$-H\"older continuous, namely,  
$$\forall z_0 \in\O, \exists M,r>0 : \r_p(z,w)\leq M\abs{z-w}^{1/2} \qquad z,w\in B(z_0,r)\subset \O.$$
\el
\bpf
Let $z_0\in \O$ and $r>0$ such that $B(z_0,r)\subset\subset \O$. Let $z,w\in B(z_0,r)$.

We use the following inequality from \cite[Proposition 4.3 (2)]{chen2022p}:
\beq \label{4.3.2}
(p-1)\abs{b-a}^2 (\abs{a}+\abs{b})^{p-2}\leq  \left[ \abs{b}^p+\abs{a}^p-\Re(\abs{b}^{p-2}\Bar{b}a+ \abs{a}^{p-2}\Bar{a}b)\right] \text{ when } p\in (1,2). 
\eeq 

Let $f_1,f_2\in \A^p(\O)$. By H\"older's inequality,
\begin{align*}
\int_\O \abs{f_1-f_2}^p &=\int_\O \abs{f_2-f_1}^p (\abs{f_1}+\abs{f_2})^\frac{p(p-2)}{2} (\abs{f_1}+\abs{f_2})^\frac{p(2-p)}{2} \\
&\leq \left[ \int_\O \abs{f_2-f_1}^2 (\abs{f_1}+\abs{f_2})^{(p-2)} \right]^{p/2}  \left[\int_\O (\abs{f_1}+\abs{f_2})^{p} \right]^{1-\frac{p}{2}}\\
&\leq \left(\frac{1}{p-1} \right)^{p/2} \left[\int_\O \abs{f_2}^p+\abs{f_1}^p-\Re(\abs{f_1}^{p-2}\Bar{f_2}f_1+ \abs{f_2}^{p-2}\Bar{f_2}f_1)\right]^{p/2} \left[\int_\O (\abs{f_1}+\abs{f_2})^{p} \right]^{1-\frac{p}{2}}.
\end{align*} Set $f_1=\frac{m_p(\z,z)}{m_p(z)}$ and $ f_2=\frac{m_p(\z,w)}{m_p(w)}$. Then
\begin{multline*}
\r_p(z,w)^p\leq \int_\O \abs{\frac{m_p(\z,z)}{m_p(z)}-\frac{m_p(\z,w)}{m_p(w)}}^p d\z \\
\hspace{-10mm}\leq \frac{2^\frac{(1-p/2)}{p}}{(p-1)^{p/2}}  \left( 1+1-\Re\left[\int_\O\abs{\frac{m_p(\z,z)}{m_p(z)}}^{p-2} \overline{\frac{m_p(\z,z)}{m_p(z)}} \frac{m_p(\z,w)}{m_p(w)} d\z +\int_\O\abs{\frac{m_p(\z,w)}{m_p(w)}}^{p-2} \overline{\frac{m_p(\z,w)}{m_p(w)}} \frac{m_p(\z,z)}{m_p(z)} d\z  \right] \right)^{p/2}.
\end{multline*}
By the reproducing formula,
$$\r_p(z,w)^p \leq \frac{2^\frac{(1-p/2)}{p}}{(p-1)^{p/2}} \Re\left[ \frac{m_p(w)-m_p(z,w)m_p(z)}{m_p(w)}+\frac{m_p(z)-m_p(w,z)m_p(w)}{m_p(z)} \right]^{p/2}.$$

Let $F(z,w)=\frac{m_p(w)-m_p(z,w)m_p(z)}{m_p(w)}$. Then $F(w,w)=0$, and from \cite{chen2022some}, we know that $m_p(\cdot,w), m_p(\cdot$) are $C^1$ functions. Thus $F(\cdot,w)$ is a $C^1$ function, hence  locally Lipschitz. Thus
$$\abs{F(z,w)-F(w,w)}=\abs{F(z,w)}\leq C_{z_0}\abs{z-w} \text{ for } z,w\in B(z_0,r) \subset\subset \O.$$
where $C_{z_0}=\max\{\abs{\frac{d}{dz}F(z,w)}: z,w\in \overline{B(z_0,r)}\}$.\\
Accordingly, $\r_p(z,w)^p \leq \frac{2^\frac{(1-p/2)}{p}}{(p-1)^{p/2}} \left[\abs{F(z,w)}+\abs{F(w,z)}\right]^{p/2} \leq C' \abs{z-w} ^{p/2}$, which implies the conclusion of the lemma. 
\epf

\section{Completeness of the \texorpdfstring{$p$}{p}-Skwarczy\'nski distance}\label{SecCompleteness}

 We'll first discuss the completeness of $p$-Skwarczy\'nski distance in the unit ball and then use some techniques of \cite{skwarczynski1980biholomorphic} to discuss the completeness of the distance in a general domain.

\bt \label{completeness}
Let $\O$ be a bounded domain in $\C^n$. Assume that for every sequence $\{w_k\}$ without an accumulation point in $\O$, 
$$\lim_{k\rightarrow \infty} \frac{m_p(z,w_k)}{m_p(w_k)} \rightarrow 0 \quad \text{ for every } z\in \O.$$ 
Then $\O$ is $p$-Skwarczy\'nski complete.
\et 

\bpf 
Let $\{w_k\}$ be a $\r_p$ Cauchy sequence. Since $\{w_k\}$ is an infinite subset of $\overline{\O}$, there exists a subsequence $w_{k_l}$ such that $w_{k_l}\xrightarrow{Euclidean} w \in \overline{\O}$.

\paragraph{Case 1.}
$w\in \O$.

Then from Lemma \ref{EucImpSkwa}, $$w_{k_l}\xrightarrow{Euclidean} w \implies w_{k_l}\xrightarrow{ \r_p} w \implies w_k\xrightarrow{\r_p} w.$$
We know that a Cauchy sequence is convergent if and only if it has a convergent subsequence, therefore $\{w_k\}$ is a convergent sequence.
\paragraph{Case 2.} $w\in \partial\O$.

By assumption, $[m_p(\cdot,w_{k_l})]$ is a Cauchy sequence in $(\P(\A^p(\O),d)$.
By completeness of $\P(\A^p(\O))$, 
there exists $f$, not identically zero, such that $[m_p(\cdot,w_{k_l})]\xrightarrow{d} [f]$, that is, there exists a sequence $\{\theta_{k_l}\}$ such that $e^{i\theta_{k_l}}\frac{m_p(\cdot,w_{k_l})}{m_p(w_{k_l})} \xrightarrow{\text{$\L_p$-norm}} f$. Consequently, these functions converge normally, and in particular pointwise.
Therefore
$$f(z)=\lim_{l\rightarrow \infty} e^{i\theta_{k_l}} \frac{m_p(z,w_{k_l})}{m_p(w_{k_l})}=0 \quad \text{ for every } z\in \O,$$
which is a contradiction. 

Combining the two cases shows that every $\r_p$ Cauchy sequence is $\r_p$ convergent.
\epf

\brem 
For the unit ball, $$m_p(\z,w)=\left[\frac{1-\abs{w}^2}{1-\<z,w\>}\right]^{4/p} \quad \text{ and } \quad m_p(w)= [\pi (1-\abs{w}^2)]^\frac{1}{p},$$ 
so the unit ball satisfies the hypothesis of Theorem \ref{completeness}. Hence \textbf{the unit ball is $p$-Skwarczy\'nski complete.}
\erem

\bl \label{MainInequality} When $p>2$, there are positive constants $c_p$ and $C_p$ such that
\beq  
\label{eq:MainInequality}
c_p \cdot {(d([m_p(\cdot,z)],[f])^p}\leq  \left[ 1- \frac{\abs{f(z)}}{K_p(z)^{1/p}}\right]\leq C_p {(d([m_p(\cdot,z)],[f])^2},
 \eeq 
 where $z\in \O$, $f\in \A^p(\O)$, and $ \norm{f}_p=1$.
\el 
The proof of Lemma \ref{MainInequality} is provided in the Appendix.\\

 Recall from \cite{skwarczynski1980biholomorphic} that if $p=2$, then
 \beq  \label{p=2}
 \frac{\abs{f(z)}}{K_p(z)^{1/p}}=\left[ 1- \frac{(d([m_p(\cdot,z)],[f])^p}{p} \right].
 \eeq

Using Lemma \ref{MainInequality}, we have the following results regarding completeness (analogous to Theorem III.6--III.13 of \cite{skwarczynski1980biholomorphic}).

\bt
A sequence $\{w_k\}$ is a $\r_p$ Cauchy sequence if and only if $\{[m_p(\cdot,w_k)]\}$ is Cauchy in $\P(\A^p)$.
\et 
\bpf
Obvious from the definition of $\r_p$.
\epf

\bt \label{Cauchy-limit}
Suppose $p>2$. A sequence $\{w_k\}$ in $\O$ is a $p$-Skwarczy\'nski Cauchy sequence if and only if there is $f$ in $\A^p(\O)$ of norm~$1$  such that  
\beq \label{Onelimit}
\lim_{k\rightarrow\infty} \frac{\abs{f(w_k)}^p}{K_p(w_k)}= 1.
\eeq
\et 

\bpf
($\impliedby$)
Assume that $$\lim_{k\rightarrow \infty} \frac{\abs{f(w_k)}}{K_p(w_k)^{1/p}} = 1.$$
By inequality \ref{eq:MainInequality},
$$c_p(d([m_p(\cdot,w_k)],[f]))^p \leq\left[ 1- \frac{\abs{f(w_k)}}{K_p(w_k)^{1/p}} \right],$$
so the assumption implies that
$$d([m_p(\cdot,w_k)],[f])\xrightarrow{k\rightarrow \infty} 0, $$
that is, $w_k$ is a $\r_p$-Skwarczy\'nski Cauchy sequence.

$(\implies)$ 
If $\{w_k\}$ is a $p$-Skwarczy\'nski Cauchy sequence, then completeness of $\P(\A^p)$ yields $f$ in $\A^p(\O)$ of norm~$1$ such that  $d([m_p(\cdot,w_k], [f])\rightarrow 0$.
By inequality \ref{eq:MainInequality}, 
$$ \left[ 1- \frac{\abs{f(z)}}{K_p(z)^{1/p}}\right]\leq C_p {(d([m_p(\cdot,z)],[f])^2}, $$
so
\begin{equation*}
    \lim_{k\rightarrow \infty} \frac{\abs{f(w_k)}}{K_p(w_k)^{1/p}}=1. \qedhere
\end{equation*}
\epf

\bt Suppose $p>2$.
Assume there exists a sequence $\{w_k\}$ with no accumulation point in $\O$ and a function $f$ in $\A^p(\O)$ of norm~$1$ such that 
\beq 
\lim_{k\rightarrow \infty} \frac{\abs{f(w_k)}^p}{K_p(w_k)}= 1.
\eeq
Then $\O$ is not $p$-Skwarczy\'nski complete.
\et

\bpf 
Suppose on the contrary that $\O$ is $p$-Skwarczy\'nski complete.
By the theorem's hypothesis and Theorem~\ref{Cauchy-limit}, the sequence $\{w_k\}$ is a $p$-Skwarczy\'nski Cauchy sequence in $\O$. 
By completeness, there is a point $w_0$ such that  $w_k\xrightarrow{\r_p} w_0$, whence  $w_k\xrightarrow{\text{Euclidean}} w_0$, contradicting the hypothesis that $\{w_k\}$ has no accumulation point.
\epf

\bt \label{ZerolimitTheorem}
Suppose $p>2$. Assume that for every sequence $\{w_k\}$ with no accumulation point in $\O$ and for every $f\in \A^p(\O)$,
\beq \label{Zerolimit}
\lim_{k\rightarrow \infty} \frac{\abs{f(w_k)}^p}{K_p(w_k)}= 0.
\eeq 
Then $\O$ is $p$-Skwarczy\'nski complete.
\et

\bpf
Suppose to the contrary that $\O$ is not $p$-Skwarczy\'nski complete, and let $\{w_k\}$ be a $p$-Skwarczy\'nski Cauchy sequence without a limit in $\O$. 
Since the two topologies agree on $\O$, the sequence $\{w_k\}$ has no accumulation point in $\O$. 
By Theorem \ref{Cauchy-limit}, there exists $f$ in $\A^p(\O)$ of norm~$1$ such that $$\lim_{k\rightarrow \infty}\frac{\abs{f(w_k)}^p}{K_p(w_k)}=1,$$ contradicting \ref{Zerolimit}. 
\epf

\bt \label{PeakingFunCriteria}
Suppose $p>2$. Assume that for each point $w\in b\O$ there exists a holomorphic peak function $h $ such that 
\begin{enumerate}[label=(\roman*)]
    \item  $\abs{h(\z)}< 1$ for every $\z$ in $\O$, and
    \item  $\lim_{\z\rightarrow w} \abs{h(\z)}=1$.
\end{enumerate}
Then $\O$ is $p$-Skwarczy\'nski complete.
\et

\bpf
Suppose to the contrary that $\O$ is not $p$-Skwarczy\'nski complete. By Theorem \ref{ZerolimitTheorem}, there exists a sequence $\{w_k\}$ without an accumulation point in $\O$ and a function $f\in \A^p(\O)$ such that 
\begin{equation} \label{eq9}
    \lim_{k\rightarrow \infty} \frac{\abs{f(w_k)}^p}{K_p(w_k)}\not \rightarrow 0.
\end{equation}  
Since $\{w_k\}$ has no accumulation point in $\O$, there exists a subsequence  $\{w_{k_l}\}$ of such that $w_{k_l}\rightarrow w \in b\O$. Thus, without loss of generality, we can assume that $w_k\rightarrow w\in b\O$ satisfying \ref{eq9}.


Given $\e>0$, notice that $\{\abs{h^kf}: k\in \Z\}$ is a sequence of functions in $\A^p(\O)$, converging pointwise to zero and is dominated by $f \in \A^p(\O)$. By dominated convergence theorem, for large enough $k_0\in \Z$, $\norm{h^{k_0}f}_p^p\leq \e$

Also, since $w_k\rightarrow w\in b\O$, for large $k_1$, $\abs{h^{k_0}(w_k)}^p\geq 1-\e$ for $k \geq k_1$

$$ (1-\e)\abs{f(w_k)}^p\leq \abs{h^{k_0}(w_k)f(w_k)}^p\leq K_p(w_k)\norm{h^{k_0}f}_p^p\leq \e K_p(w_k) \text{ for } k\geq k_1$$

Thus $$\frac{\abs{f(w_k)}^p}{K_p(w_k)}\leq \frac{\e}{1-\e} \text{ for } k\geq k_1 \implies \lim_{k\rightarrow \infty } \frac{\abs{f(w_k)}^p}{K_p(w_k)}\rightarrow 0$$
contradicting our assumption
\epf

\bcor
Let $G$ be a domain in $\C^n$ and $\O$ be a connected analytic polyhedron defined by its frame $\{h_i\}_{i=1}^k \subset \mathcal{O}(G)$ $$\O=\{z\in G: \abs{h_i(z)}<1 ; i=1,\dots,k\} .$$ Then $\O$ is $p$-Skwarczy\'nski complete, when $p>2$.
\ecor

\bpf
For every $w_0\in b\O$, there is at least one $i=i_0$ such that $\abs{h_{i_0}(w_0)}=1$. Thus $\O$ satisfies the hypothesis of Theorem \ref{PeakingFunCriteria} and hence $\O$ is $p$-Skwarczy\'nski complete for $p>2$.
\epf

\brem
The following domains are known to satisfy the hypothesis of Theorem \ref{PeakingFunCriteria} :
\begin{enumerate}
    \item Smooth bounded strictly pseudoconvex domains \cite[Prop 2.1]{noell2008peak};
    \item $h$-extendible (semiregular) domains \cite[Theorem A]{yu1994peak};
    \item Bounded pseudoconvex domains in $\C^2$ with a real analytic boundary \cite[Theorem 3.1]{bedford1978construction}. 
\end{enumerate}
Thus, for any domain $\O$ of the above type, $\O$ is $p$-Skwarczy\'nski complete for $p>2.$\erem

\bt \label{DenseBoundedFn}
Suppose $p>2$. Assume that for every point $w\in b\O$,  
\begin{enumerate}[label=(\roman*)]
    \item $\lim_{\z \rightarrow w} K_p(\z)\rightarrow \infty$, and
    \item $\mathcal{O}(\O)\cap C(\overline{\O})$ is dense  in $\A^p(\O)$. 
\end{enumerate}
Then $\O$ is $p$-Skwarczy\'nski complete.
\et
\bpf
We verify the hypothesis of Theorem \ref{ZerolimitTheorem}. 
Let $f\in \A^p(\O)$ be any holomorphic function and let $\{z_k\}$ be any sequence in $\O$ without an accumulation point in $\O$. Then we can find a sequence of holomorphic functions $g_k \in \mathcal{O}(\O)\cap C(\overline{\O})$ such that  $g_k\xrightarrow{\norm{.}_p} f$. Hence we have $$\lim_{k\rightarrow \infty} \frac{\abs{f(z_k)}^p}{K_p(z_k)} =\lim_{k\rightarrow \infty} \frac{\abs{g_k(z_k)}^p}{K_p(z_k)}=0 $$
 proving $\O$ is $p$-Skwarczy\'nski complete.

\epf



\section{Invariance of the \texorpdfstring{$p$}{p}-Skwarczy\'nski distance}\label{SecInvariance}

\subsection{\texorpdfstring{$\A^p$}{Ap} preserving biholomorphisms}

Let  $F: \O_1\rightarrow \O_2$ be a biholomorphism. Then $F$ induces an isometric isomorphism on $\A^2$ in a natural way, namely 
 $$ \hspace{-23mm}F^{\#}:\A^2(\O_2) \rightarrow \A^2(\O_1) $$
 $$f \mapsto f\circ F \cdot J_F$$ where $J_F=det(F')$.\\
 By the change of variables formula, $$\hspace{-25mm}\int_{\O_1}  \abs{f\circ F}^2 \abs{J_F}^2 =\int_{\O_2} \abs{f}^2$$  thereby defining an isometric isomorphism. However in case of $\A^p$, not all biholomorphisms have a `natural' extension onto the $\A^p$ space. 
 
 \bedef
 Let $F:\O_1 \rightarrow \O_2$ be a biholomorphism between two bounded domains $\O_1,\O_2$. We say that $F$ is an \textbf{\texorpdfstring{$\A^p$}{Ap} preserving biholomorphism} if the map $F^\#$ extends naturally onto $\A^p$ space defining an isometric isomorphism i.e.
 $$\hspace{-25mm}F^{\#}:\A^p(\O_2) \rightarrow \A^p(\O_1) $$
 $$f \mapsto f\circ F \cdot J_F^{2/p}$$ 
 $$\norm{f}_{\A^p(\O_2)}= \norm{f\circ F \cdot J_F^{2/p}}_{\A^p(\O_1)}$$
 \endef 
 
 \brem
\begin{itemize}
  \item The above definition is equivalent to saying the $p/2^{th}$ root of the function $J_F$ is a holomorphic function on $\O_1$.
  \item Every biholomorphism is $\A^2$ preserving and $\A^p$ preserving for $p$ such that $2/p\in \N$.
  \item Every biholomorphism on a simply connected domain is an $\A^p$ preserving biholomorphism.
  \item If $F:\O_1\rightarrow \O_2$ is an $\A^p$ preserving biholomorphism, then $F^{-1}:\O_1\rightarrow \O_2$ is also an $\A^p$ preserving biholomorphism.
\end{itemize}
 \erem 
 
\bt
Let $p\geq1$ and $F:\O_1\rightarrow \O_2$ be an $\A^p$ preserving biholomorphism between two bounded domains $\O_1,\O_2$. Then $$\r_{p,\O_1}(z,w)=\r_{p,\O_2}(F(z),F(w)) \text{ for } z,w\in \O_1.$$
In other words, the $p$-Skwarczy\'nski distance is invariant under $\A^p$ preserving biholomorphisms.
\et 
\bpf
Let $F:\O_1\rightarrow \O_2$ be an $\A^p$ preserving biholomorphsim. 
Then from \cite{chen2022p}, we have 
\begin{align*}
m_{p,\O_1}(w)&=m_{p,\O_2}(F(w))\cdot \abs{J_F(w)}^{-2/p}.\\
m_{p,\O_1}(\z,w)&=m_{p,\O_2}(F(\z),F(w))\cdot J_F(w)^{-2/p}\cdot J_F(\z)^{2/p}.
\end{align*}

By definition,
\begin{align*}
&\hspace{-5mm}\r_{p,\O_1}^p(z,w)=\min_{t_1,t_2\in [0,2\pi]} \norm{e^{it_1}\frac{m_{p,\O_1}(\cdot,w)}{m_{p,\O_1}(w)}- e^{it_2}\frac{m_{p,\O_1}(\cdot,z)}{m_{p,\O_1}(z)}}^p_{p,\O_1}\\
&\hspace{-7mm}=\min_{t_1,t_2\in [0,2\pi]} \int_{\O_1} \abs{e^{it_1}\frac{m_{p,\O_2}(F(\z),F(w))\cdot J_F(w)^{-2/p}\cdot J_F(\z)^{2/p}} {m_{p,\O_2}(F(w))\cdot \abs{J_F(w)}^{-2/p}}- e^{it_2}\frac{m_{p,\O_2}(F(\z),F(z))\cdot J_F(z)^{-2/p}\cdot J_F(\z)^{2/p}}{m_{p,\O_2}(F(z))\cdot \abs{J_F(z)}^{-2/p}}}^pd\z\\
&\hspace{-7mm}= \min_{t_1',t_2'\in [0,2\p]} \int_{\O_2} \abs{ e^{it_1'}\frac{m_{p,\O_2}(\z',F(w))}{m_{p,\O_2}(F(w))}- e^{it_2'}\frac{m_{p,\O_2}(\z',F(z))}{m_{p,\O_2}(F(z))}}^p d\z'= \r^p_{p,\O_2}(F(z),F(w))
\end{align*}

Thus the $p$-Skwarczy\'nski distance is invariant under $\A^p$ preserving biholomorphisms.
\epf

\section{Continuity of the \texorpdfstring{$p$}{p}-Skwarczy\'nski distance \w.r.t. \texorpdfstring{$p$}{p}} \label{SecContinuity}
We will use the following result from \cite[Theorem 2.2]{xiangrevisit}. We will prove it here for the sake of completeness.
\bl \label{BrezisExtended}
Let $\{f_k\} \subset \L^q(\O)$ for  some $1\leq q<\infty $. Assume that  \begin{itemize}
    \item[i.] $f_k\rightarrow f$ pointwise almost everywhere on $\O$
    \item[ii.] the sequence $\{f_k\}$ is uniformly bounded in $\L^q$ i.e. $\norm{f_k}_q\leq M$ for all $k$ and some $M>0$.
\end{itemize}
If the volume of $\O$, $\abs{\O}$,  is finite, then $f\in \L^q(\O)$ and $f_k\rightarrow f$ in $\L^p(\O)$ for $0<p<q.$
\el
\bpf
By Fatou's Lemma, 
$$\norm{f}_q \leq \liminf_{k\rightarrow \infty} \norm{f_k}_q \leq M \implies f\in \L^q(\O).$$
Let $E\subset \O $ be a measurable set and $0<p<q$. Applying H\"older's inequality and triangle inequality 
$$ \int_E \abs{f_k-f}^p \leq \norm{f_k-f}_q^p \abs{E}^{\frac{q-p}{q}} \leq (M+ \norm{f}_q)^p \abs{E}^{\frac{q-p}{q}} $$
Thus $\abs{f_k-f}^p \rightarrow 0$ a.e. and is uniformly integrable over $\O$. By Vitali's convergence theorem $f_k\rightarrow f$ in $\L^p(\O)$.
\epf 

A bounded domain $\O$ is said to be hyperconvex if there exists a negative continuous plurisubharmonic function $r$, such that $\{r<c\}\subset \subset \O$ for all $c<0$. Further assume that the above function $r$ satisfies a growth condition $-r\leq C \d^\a$,  for some $\a, C>0$ and $\d$ denotes the distance function to the boundary. Let $\a(\O)$ be the supremum of all such $\a$. The hyperconvexity index $\a(\O)$ is studied in detail in  \cite{chen2017bergman}

\bt 
Let $\O$ be a bounded hyperconvex domain with $\a(\O)>0$. 
\begin{itemize}
    \item[a.]  Then for $p\in (1,2]$ $$\lim_{q\rightarrow p^-} \r_q(z,w)= \r_p(z,w) $$
    \item[b.] Let $p \in [1,2)$. Additionally if $\O$ is such that $A^{p'}(\O)$ is dense in $A^p(\O)$ for some $p'>p$, then $$\lim_{s\rightarrow p^+} \r_s(z,w)=\r_p(z,w)$$ 
\end{itemize}
\et 

\bpf
Fix $p\in (1,2]$. Suppose that $\O$ is a bounded hyperconvex domain with $\a(\O)>0$. Then from \cite[Theorem 1.4]{chen2022some}
\begin{equation} \label{hyperconvex domain}
K_p(\cdot,z)\in \L^q(\O) \text{ for } q<\frac{2pn}{2n-\a(\O)}
\end{equation}
and from \cite[Theorem 6.5]{chen2022p} 
$$\lim_{s\rightarrow p^-} m_s(\z,w)= m_p(\z,w) \text{ for } z,w\in \O $$
Let $p_n\nearrow p$ be any sequence and let $$f_n^\theta(\z)=e^{i \theta} \frac{m_{p_n}(\z,z)}{m_{p_n}(z)}- \frac{m_{p_n}(\z,w)}{m_{p_n}(w)}$$
$$f^\theta(\z)=e^{i \theta} \frac{m_{p}(\z,z)}{m_{p}(z)}- \frac{m_{p}(\z,w)}{m_{p}(w)}$$
Choose $\theta_n$ such that $\r_{p_n}(z,w)=\norm{f_n^{\theta_n}}_{p_n }$ and $\t_0$ such that $\r_p(z,w)=\norm{f^{\t_0}}_p$. \\
(We know that $\r$ is a continuous function if for every sequence $x_n \rightarrow x$, we have a subsequence of $\r(x_n)$ which converges to $\r(x)$.) Therefore, without loss of any generality, we can assume that $\theta_n \rightarrow \theta_0 $ for some $\theta_0 \in [0,2\pi]$.\\\\
By \cite[Theorem 1.4]{chen2022some}, for $p_n$ sufficiently close to $p$, $f_n \in L^{q}(\O)$ (for some $q>p$) and by \cite[Theorem 6.5]{chen2022p} $ f_n^{\theta_n} \xrightarrow{pointwise} f^{\theta_0}$ and $ f_n^{\t_0} \xrightarrow{pointwise} f^{\t_0} $.
By Lemma \ref{BrezisExtended}, for every $\e>0$, we can find $N_0\in \N$ such that 
\begin{align*}
    \norm{f^{\theta_0}-f_n^{\theta_n}}_{p} \leq \e/2 &\text{ for } n\geq N_0\\
    \norm{f_n^{\t_0}}_{p} \leq \norm{f^{\t_0}}_{p} +\e/2 &\text{ for } n\geq N_0
\end{align*}
Let $p'<p$. Choose $p'<p_k<p$  and $k\geq N_0$. Using the triangle inequality and  H\"older's inequality, 
\begin{align*}
    \norm{f^{\theta_0}}_{p'} &\leq \norm{f^{\theta_0}-f_k^{\theta_k}}_{p'}+\norm{f_k^{\theta_k}}_{p'}\\
    & \leq \norm{f^{\theta_0}-f_k^{\theta_k}}_{p}\cdot\abs{\O}^{\frac{1}{p'}-\frac{1}{p}}+ \norm{f_k^{\theta_k}}_{p_k}\cdot\abs{\O}^{\frac{1}{p'}-\frac{1}{p_k}} \\
    & \leq \norm{f^{\theta_0}-f_k^{\theta_k}}_{p}\cdot\abs{\O}^{\frac{1}{p'}-\frac{1}{p}}+ \norm{f_k^{\t_0}}_{p_k}\cdot\abs{\O}^{\frac{1}{p'}-\frac{1}{p_k}}\\
    & \leq \norm{f^{\theta_0}-f_k^{\theta_k}}_{p}\cdot\abs{\O}^{\frac{1}{p'}-\frac{1}{p}}+ \norm{f_k^{\t_0}}_{p}\cdot\abs{\O}^{\frac{1}{p'}-\frac{1}{p}}\\
    & \leq (\e/2+\e/2+\norm{f^{\t_0}}_{p}) \cdot \abs{\O}^{\frac{1}{p'}-\frac{1}{p}}
\end{align*}
By continuity of norm $$\norm{ f^{\theta_0}}_{p}\leq \norm{f^{\t_0}}_{p}+\e \text{ and } \norm{f^{\t_0}}_p\leq \norm{f^{\theta_0}}_p$$ 
thereby proving, $$\lim_{q\rightarrow p^-} \r_{q}(z,w)=\norm{f^{\theta_0}}_p=\norm{f^{\t_0}}_p=\r_p(z,w)$$

Additionally, assume that $A^{q}(\O)$ lies dense in $A^p(\O)$ for some $q>p$. Then by \cite[Theorem 6.5]{chen2022p}
\begin{equation}
    \lim_{ s\rightarrow p^+}m_s(\z,w)= m_p(\z,w)
\end{equation}
Let $1\leq p< 2$, $p_n\searrow p$ $(p_1<2)$ be any sequence and 
$$g_n^{\theta}(\z)=e^{i \theta} \frac{m_{p_n}(\z,z)}{m_{p_n}(z)}- \frac{m_{p_n}(\z,w)}{m_{p_n}(w)}$$
$$g^\theta(\z)=e^{i \theta} \frac{m_{p}(\z,z)}{m_{p}(z)}- \frac{m_{p}(\z,w)}{m_{p}(w)}$$

Choose $\theta_n$ such that $\r_{p_n}(z,w)=\norm{g_n^{\theta_n}}_{p_n }$ and $\t_0$ such that $\r_p(z,w)=\norm{g^{\t_0}}_p$. As above, without loss of any generality, we can assume that $\theta_n \rightarrow \theta_0 $ for some $\theta_0 \in [0,2\pi]$.\\\\
According to \cite[Theorem 1.4]{chen2022some}, $g_n^{\theta}\in L^{p_1}(\O)$ and by \cite[Theorem 6.5]{chen2022p} $g_n^{\theta_n}\xrightarrow{pointwise}g^{\theta_0}$ and $g_n^{\t_0}\xrightarrow{pointwise}g^{\t_0}$.
Let $p'>p$. By Lemma \ref{BrezisExtended}, for every $\e>0$, we can find $N_0\in \N$ such that 
\begin{align*}
    \norm{g^{\theta_0}-g_n^{\theta_n}}_{p'} \leq \e/2 &\text{ for } n\geq N_0\\
    \norm{g_n^{\t_0}}_{p'} \leq \norm{g^{\t_0}}_{p'} +\e/2 &\text{ for } n\geq N_0
\end{align*}
Choose $p<p_k<p'$, $k\geq N_0$. Using triangle inequality and  H\"older's inequality, 
\begin{align*}
    \norm{g^{\theta_0}}_{p} &\leq \norm{g^{\theta_0}-g_k^{\theta_k}}_{p}+\norm{g_k^{\theta_k}}_{p}\\
    & \leq \norm{g^{\theta_0}-g_k^{\theta_k}}_{p'}\cdot\abs{\O}^{\frac{1}{p}-\frac{1}{p'}}+ \norm{g_k^{\theta_k}}_{p_k}\cdot\abs{\O}^{\frac{1}{p}-\frac{1}{p_k}} \\
    & \leq \norm{g^{\theta_0}-g_k^{\theta_k}}_{p'}\cdot\abs{\O}^{\frac{1}{p}-\frac{1}{p'}}+ \norm{g_k^{\t_0}}_{p_k}\cdot\abs{\O}^{\frac{1}{p}-\frac{1}{p_k}}\\
    & \leq \norm{g^{\theta_0}-g_k^{\theta_k}}_{p'}\cdot\abs{\O}^{\frac{1}{p}-\frac{1}{p'}}+ \norm{g_k^{\t_0}}_{p'}\cdot\abs{\O}^{\frac{1}{p'}-\frac{1}{p}}\\
    & \leq (\e/2+\e/2+\norm{g^{\t_0}}_{p'}) \cdot \abs{\O}^{\frac{1}{p'}-\frac{1}{p}}
\end{align*}
By continuity of norm $$\norm{ g^{\theta_0}}_{p}\leq \norm{g^{\t_0}}_{p}+\e \text{ and } \norm{g^{\t_0}}_p\leq \norm{g^{\theta_0}}_p$$ 
which proves, $$\lim_{s\rightarrow p^+} \r_{s}(z,w)=\norm{g^{\theta_0}}_p=\norm{g^{\t_0}}_p=\r_p(z,w)$$
\epf

\section{Product Domains} \label{SecProductDomain}

\subsection{\texorpdfstring{$p$}{p}-Skwarczy\'nski distance on the product domain}

\bl

Suppose that $\O_1 \subset \C^{n_1}, \O_2\subset \C^{n_2}$ are two bounded domains. Let $\O=\O_1 \times \O_2$, $z=(z_1,z_2)$ , $ w=(w_1,w_2)$. Then $$\r_{p,\O} (z,w)\leq \r_{p,\O_1}(z_1,w_1)+\r_{p,\O_2}(z_2,w_2)$$ 

\el 

\bpf 
By definition and the product rule from \cite[Proposition 2.8]{chen2022p}
$$ \hspace{-9mm}\r_{p,\O}(z,w)=\min_{\theta\in [0,2\pi]}\norm{e^{i\theta}\frac{m_{p,\O}(\cdot,z)}{m_{p,\O}(z)}-\frac{m_{p,\O}(\cdot,w)}{m_{p,\O}(w)}}_{p,\O}$$
$$ \hspace{55mm}=\min_{\theta\in [0,2\pi]}\norm{e^{i\theta}\frac{m_{p,\O_1}(\cdot,z_1)m_{p,\O_2}(\cdot,z_2)}{m_{p,\O_1}(z_1)m_{p,\O_2}(\cdot,z_2)}-\frac{m_{p,\O_1}(\cdot,w_1)m_{p,\O_2}(\cdot,w_2)}{m_{p,\O_1}(w_1)m_{p,\O_2}(w_2)}}_{p,\O}$$

For $i=1,2$,    $$\text{set } f_i(\z_i)=e^{i\theta_i}\frac{m_{p,\O_i}(\z_i,z_i)}{m_{p,\O_i}(z_i)} ; \hspace{3mm} g_i(\z_i)=\frac{m_{p,\O_i}(\z_i,w_i)}{m_{p,\O_i}(w_i)} \in \A^p(\O_i)$$ where $\theta_i$ is chosen such that $\r_{p,\O_i}(z_i,w_i)=\norm{f_i-g_i}_{p,\O_i}$.\\
Then $$\r_{p,\O}(z,w)\leq \norm{f_1f_2-g_1g_2}_{p,\O}\leq  \norm{f_1f_2-g_1f_2}_{p,\O}+\norm{g_1f_2-g_1g_2}_{p,\O}.$$
Consider $$\norm{f_1f_2-g_1f_2}_{p,\O}=\left(\int_{\O_1\times\O_2} \abs{f_1(\z_1)-g_1(\z_1)}^p \abs{f_2(\z_2)}^p d\z_1d\z_2\right)^{1/p}$$
$$\hspace{60mm}=\left(\int_{\O_1} \abs{f_1(\z_1)-g_1(\z_1)}^p  d\z_1\right)^{1/p} \left(\int_{\O_2}  \abs{f_2(\z_2)}^p d\z_2\right)^{1/p}=\r_{p,\O_1}(z_1,w_1)$$
Similarly $\norm{g_1f_2-g_1g_2}_{p,\O}=\r_{p,\O_2}(z_2,w_2)$, so
$$\r_{p,\O} (z,w)\leq \r_{p,\O_1}(z_1,w_1)+\r_{p,\O_2}(z_2,w_2)$$ 
\epf

\subsection{\texorpdfstring{$p$}{p}-Bergman metric on the product domain}

In \cite{chen2022p}, the $p$-Bergman metric was defined as follows for a vector field~$X$.
$$ B_p(z_0;X):= K_p(z_0)^{-\frac{1}{p}} \sup\{\,\abs{Xf(z_0)}: f\in \A^p, \; f(z_0)=0, \; \norm{f}_p=1\,\}.$$ 

The problem was posed to find the $p$-Bergman metric on product of two domains.
Here, we provide a partial solution to this problem.

\bl 
Suppose that $\O_1 \subset \C^{k_1}$ and $\O_2\subset \C^{k_2}$ are bounded domains. Let $\O=\O_1 \times \O_2$, $z=(z_1,z_2)$, $ X=(X_1,X_2)$. Then $$B_{p,\O} (z;X)\geq \max_{i=1,2} B_{p,\O_i}(z_i;X_i).$$ 
\el 

\bpf 
Suppose $h\in \A^p(\O_1)$, $\norm{h}_{p,\O_1}=1$,  $f\in \A^p(\O_2)$, $f(z_2)=0$, and $\norm{f}_{p,\O_2}=1$.
If $g(\z_1,\z_2):=h(\z_1)\cdot f(\z_2) $, then $g\in \A^p(\O)$, and $\norm{g}_{p,\O}=\norm{h}_{p,\O_1}\cdot \norm{f}_{p,\O_2}=1$.
Now
$$B_{p,\O}(z;X)\geq K_{p,\O}^{-\frac{1}{p}}(z)  \abs{X(g)(z)}=K_{p,\O_1}(z_1)^{-\frac{1}{p}}\abs{h(z_1)} \cdot K_{p,\O_2}(z_2)^{-\frac{1}{p}}\abs{ X_2(f)(z_2)}.$$
Taking the supremum over all functions $f$ and $h$ as above, we get $$\sup_{h}K_{p,\O_1}(z_1)^{-\frac{1}{p}}\abs{h(z_1)}=1 \text{ and } \sup_{f} K_{p,\O_2}(z_2)^{-\frac{1}{p}} \abs{X_2(f)(z_2)}= B_{p,\O_2}(z_2;X_2). $$  
By symmetry,
\begin{equation*}B_{p,\O} (z;X)\geq \max_{i=1,2} B_{p,\O_i}(z_i;X_i). \qedhere
\end{equation*}
\epf

\section{Application} \label{SecApplications}

From \cite{chen2022p}, we have $\abs{m_p(z,w)}\leq m_p(w)/m_p(z)$, and equality holds if and only if $z=w$.
This follows from the reproducing formula and H\"older's inequality. Indeed,
$$\abs{m_p(z,w)}=\abs{m_p(z)^{-p}\int_\O\abs{m_p(\z,z)}^{p-2}\overline{m_p(\z,z)}m_p(\z,w) d\z  }\leq m_p(w)/m_p(z).$$

We can find a stronger bound using the $p$-Skwarczy\'nski distance.

\bl 
 If $p> 2$,  then \begin{equation}\label{inequality}
 \abs{m_p(z,w)}\leq \frac{m_p(w)}{m_p(z)} \left[ 1- \frac{\r_p(z,w)^p}{p\cdot 4^{p+3}} \right].\end{equation}
\el 
Equation \eqref{inequality} shows that 
$\abs{m_p(z,w)}=m_p(w)/m_p(z)$ if and only if $\r_p(z,w)=0$ (equivalently $z=w$). 

\bpf
From \cite[Proposition 4.3 (3)]{chen2022p}, we have $$\abs{b}^p \geq \abs{a}^p +p\Re(\abs{a}^{p-2}\bar{a}(b-a))+\frac{1}{4^{p+3}} \abs{b-a}^p \text{ when }p>2. $$
Set $a=m_p(\z,z)/m_p(z)$, and
$b=e^{i\theta }m_p(\z,w)/m_p(w)$, where $\theta$ will be specified below.  Integrating the above inequality shows that
\begin{multline*}
    1\geq 1 + p \Re\left\{ \int_\O m_p(z)^{-p+1}\abs{m_p(\z,z)}^{p-2} \overline{m_p(\z,z)} \left[\frac{e^{i\theta }m_p(\z,w)}{m_p(w)}- \frac{m_p(\z,z)}{m_p(z)} \right] d\z \right\} \\ + \frac{1}{4^{p+3}}\int_\O \abs{ \frac{e^{i\theta }m_p(\z,w)}{m_p(w)}- \frac{m_p(\z,z)}{m_p(z)} }^p d\z.\end{multline*}
Applying the reproducing property shows that
$$ p m_p(z)^{} \Re\left\{ \frac{e^{i\theta }m_p(z,w)}{m_p(w)}- \frac{m_p(z,z)}{m_p(z)}\right\} +\frac{\r_p(z,w)^p}{4^{p+3}}\leq 0.$$
Now choose $\theta$ such that $e^{i\theta} m_p(z,w)=\abs{m_p(z,w)}$.
Then, $$\frac{\abs{m_p(z,w)} m_p(z)}{m_p(w)} -1 \leq- \frac{\r_p(z,w)^p}{p\cdot 4^{p+3}}$$
which is equivalent to \eqref{inequality}.
\epf 

\section{Appendix}

 We will now look at the proof of Lemma \ref{MainInequality} by partially following the steps in \textit{Appendix} of \cite{chen2022p}.
\bpf
 Let $a,b\in \C$, $p\geq 1$. Define $$\eta(t)=\abs{a+t(b-a)}^2 ; \kappa(t)=\eta(t)^{p/2} =\abs{a+t(b-a)}^p.$$

$$\k'(t)=\frac{p}{2} \eta(t)^{p/2-1}\eta'(t) = p\cdot\abs{a+t(b-a)}^{p-2} Re(\bar{a}(b-a)+t(\abs{a-b}^2)).$$
Using  $(Re\{\bar{a}(b-a)\}+ t\abs{b-a}^2)^2 +(Im(\bar{a}b))^2=\abs{b-a}^2\abs{a+t(b-a)}^2$
$$\k''(t)=p\abs{a+t(b-a)}^{p-4} \left[(Im\{\bar{a}b\})^2+(p-1) \cdot (Re\{\bar{a}(b-a)\}+ t\abs{b-a}^2)^2 \right]$$
which implies, 
$$ p\min\{1,(p-1)\}\abs{a+t(b-a)}^{p-2}\abs{b-a}^2 \leq \k''(t)\leq p\max\{1,(p-1)\}\abs{a+t(b-a)}^{p-2}\abs{b-a}^2  $$

Using integration by parts, we have 

$$\k(1)=\k(0)+\k'(0)+\int_0^1(1-t)k''(t)$$

Applying the upper and lower limits on $\k''(t)$, we get 

\beq \label{MI1}
\abs{b}^p \geq \abs{a}^p +pRe(\abs{a}^{p-2}\bar{a}(b-a))+ p\min\{1, (p-1)\}\int_0^1(1-t) \abs{a+t(b-a)}^{p-2}\abs{b-a}^2dt  
\eeq 

\beq \label{MI2}
\abs{b}^p \leq \abs{a}^p +pRe(\abs{a}^{p-2}\bar{a}(b-a))+  p\max\{1, (p-1)\}\int_0^1 (1-t) \abs{a+t(b-a)}^{p-2}\abs{b-a}^2dt
\eeq

Let $p>2$, using \eqref{MI1} $$\abs{b}^p \geq \abs{a}^p +pRe(\abs{a}^{p-2}\bar{a}(b-a))+ p\cdot \abs{b-a}^2\cdot I, $$
 where $I=\int_0^1 (1-t)\abs{a+t(b-a)}^{p-2}dt$.

We will now compute a lower bound on $I$
$$I\geq \int_0^1 \abs{\abs{a}-t\abs{b-a}}^{p-2}dt.$$

If $a\geq \abs{b-a}/2$, then
$$I\geq \abs{b-a}^{p-2} \int_0^{1/4} (1-t)(1/2-t)^{p-2}dt. $$

If $a\leq \abs{b-a}/2$, then
$$I\geq \abs{b-a}^{p-2} \int_{3/4}^{1} (1-t)(t-1/2)^{p-2} dt. $$

Thus there exists $c>0 $ such that $I\geq c \abs{b-a}^{p-2}$. Hence we have

$$\abs{b}^p \geq \abs{a}^p +pRe(\abs{a}^{p-2}\bar{a}(b-a))+c\abs{b-a}^p. $$

Let $f\in \A^p, \norm{f}_p=1$ $b=e^{i\theta }f(\z), a=m_p(\z,z)/m_p(z)$ and integrate the above inequality

$$1\geq \frac{m_p(z)^p}{m_p(z)^p}+ p Re\left\{ \int_\O m_p(z)^{-p+1}\abs{m_p(\z,z)}^{p-2} \overline{m_p(\z,z)} \left[e^{i\theta}f(\z)- \frac{m_p(\z,z)}{m_p(z)} \right] d\z \right\} $$
$$ + c\int_\O \abs{ \frac{e^{i\theta }m_p(\z,w)}{m_p(w)}- \frac{m_p(\z,z)}{m_p(z)} }^p d\z$$

Using the reproducing property,
$$ p m_p(z)^{} Re\left\{ e^{i\theta}f(z)- \frac{1}{m_p(z)}\right\} +c\cdot d([m_p(\cdot,z)],[f])^p\leq 0. $$

Choose $\theta$ such that $e^{i\theta} f(z)=\abs{f(z)}$

 \beq  \label{Inequality}
 \abs{f(z)}m_p(z)\leq\left[ 1- \frac{(d([m_p(\cdot,z)],[f])^p}{c_p} \right]
 \eeq

Let $p>2$. Using \eqref{MI2}
 
 $$\abs{b}^p \leq \abs{a}^p+p\cdot Re\left[ \abs{a}^{p-2}\bar{a}(b-a)\right]+p(p-1)\int_{0}^1 \abs{b-a}^2 [(1-t)]\abs{a+t(b-a)}^{p-2}dt $$

 Let $f\in \A^p$, $\norm{f}_p=1$, $b=e^{i\theta }f(\z)$, $a=m_p(\z,z)/m_p(z)$ and integrate the above inequality
 
 $$\int_{\O} f(\z)^p \leq \int_\O \abs{\frac{m_p(\z,z)}{m_p(z)}}^p +p Re\left\{ \int_\O m_p(z)^{-p+1}\abs{m_p(\z,z)}^{p-2} \overline{m_p(\z,z)} \left[e^{i\theta}f(\z)- \frac{m_p(\z,z)}{m_p(z)} \right] d\z \right\}$$  $$ +  p(p-1) \int_\O \int_{0}^1 (1-t) \abs{e^{i\theta}f(\z)-\frac{m_p(\z,z)}{m_p(z)}}^2 \abs{\frac{m_p(\z,z)}{m_p(z)}+t\left[e^{i\theta}f(\z)-\frac{m_p(\z,z)}{m_p(z)}\right]}^{p-2}dt d\z $$

 Consider $$I_1=\int_\O \int_{0}^1 [(1-t)] \abs{e^{i\theta}f(\z)-\frac{m_p(\z,z)}{m_p(z)}}^2 \abs{\frac{m_p(\z,z)}{m_p(z)}+t\left[e^{i\theta}f(\z)-\frac{m_p(\z,z)}{m_p(z)}\right]}^{p-2}dt d\z$$ 
 
 Using Fubini's theorem and H\"older's inequality (with $p'=p/2, q'=p/(p-2)$)
 
 $$I_1\leq \int_{0}^1 [(1-t)] \left[ \int_\O\abs{e^{i\theta}f(\z) -\frac{m_p(\z,z)}{m_p(z)}}^p d\z \right]^{2/p} \left[ \int_\O \abs{\frac{m_p(\z,z)}{m_p(z)}+t\left[e^{i\theta}f(\z)-\frac{m_p(\z,z)}{m_p(z)}\right]}^{p} d\z \right]^{(p-2)/p}dt$$
 
 $$\hspace{-87mm}\leq \frac{3^{(p-2)/p}}{2} \left[ \int_\O\abs{e^{i\theta}f(\z) -\frac{m_p(\z,z)}{m_p(z)}}^p d\z \right]^{2/p} $$
 Thus using the reproducing formula, $$p\cdot Re\left[ m_p(z) e^{i\theta} f(z)-1\right] + C_p \left[ \int_\O\abs{e^{i\theta}f(\z) -\frac{m_p(\z,z)}{m_p(z)}}^p d\z \right]^{2/p} \geq 0 $$ for all $\theta$ and for some $C_p>0$.\\\\ Choose $\theta$ such that $\int_\O\abs{e^{i\theta}f(\z) -\frac{m_p(\z,z)}{m_p(z)}}^p d\z =(d([m_p(\cdot,z)],[f])^p$. Then,
 
 $$C' \cdot d([m_p(\cdot,z)],[f])^2 \geq 1- m_p(z) Re\{e^{i\theta}f(z)\}\geq 1-m_p(z)\abs{f(z)}$$

 \epf

\subsection*{Acknowledgements}
{\fontsize{11.5}{10}\selectfont
The author would like to thank Prof.\  Harold Boas for giving valuable advice and feedback during the preparation of this note. He would also like to thank Tanuj Gupta, Siddharth Sabharwal and John Treuer for useful conversations.}

\bibliographystyle{unsrt}  
\bibliography{references}

\fontsize{11}{9}\selectfont

\vspace{0.5cm}

\noindent bhatshreedhar33@tamu.edu;

 \vspace{0.2 cm}

\noindent Department of Mathematics, Texas A\&M University, College Station, TX 77843-3368, USA

\end{document}